\documentclass[12pt]{amsart}
\usepackage{amssymb}
\begin{document}

\def\dbl{[\hskip -1pt[}
\def\dbr{]\hskip -1pt]}
\newcommand{ \al}{\alpha}
\newcommand{\ab}[1]{\vert z\vert^{#1}}
\title
{
Nonlinear  CR automorphisms of Levi degenerate hypersurfaces and a
new gap phenomenon.
}
\author{ martin kolar and francine meylan }

\address{M. Kolar: Department of Mathematics and Statistics, Masaryk University,
Kotlarska 2,  611 ~37
 Brno, Czech Republic}

\email{mkolar@math.muni.cz}
\thanks{The first author was supported by the project CZ.1.07/2.3.00/20.0003
of the Operational Programme Education for Competitiveness of the Ministry
of Education, Youth and Sports of the Czech Republic.}

\address{ F. Meylan: Department of Mathematics,
University of Fribourg, CH 1700 Perolles, Fribourg}

\email{francine.meylan@unifr.ch}
\thanks{The second author was supported by Swiss NSF Grant 2100-063464.00/1 }

\begin{abstract}
We give a complete classification of polynomial models for smooth
real hypersurfaces of finite Catlin multitype in $\mathbb C^3$,
which admit nonlinear
infinitesimal CR automorphisms.
As a consequence, we obtain a sharp 1-jet determination result for
any smooth hypersurface with such model.
The results also prove
a conjecture of the first author about the origin of such
nonlinear automorphisms (AIM list of problems, 2010).
 As another  consequence, we
describe all possible dimensions of the Lie algebra of
infinitesimal CR automorphisms, which leads to
a new
"secondary" gap phenomenon.
\end{abstract}


\def\Label#1{\label{#1}}
\def\1#1{\ov{#1}}
\def\2#1{\widetilde{#1}}
\def\6#1{\mathcal{#1}}
\def\4#1{\mathbb{#1}}
\def\3#1{\widehat{#1}}
\def\K{{\4K}}
\def\LL{{\4L}}
\def\H{{\4H}}
\def\C{{\4C}}
\def\R{{\4R}}
\def \MM{{\4M}}
\def\La{\Lambda}
\def\la{\lambda}
\def\Re{{\sf Re}\,}
\def\Im{{\sf Im}\,}

\numberwithin{equation}{section}
\def\s{s}
\def\k{\kappa}
\def\ov{\overline}
\def\span{\text{\rm span}}
\def\ad{\text{\rm ad }}
\def\tr{\text{\rm tr}}
\def\xo {{x_0}}
\def\Rk{\text{\rm Rk\,}}
\def\sg{\sigma}
\def \emxy{E_{(M,M')}(X,Y)}
\def \semxy{\scrE_{(M,M')}(X,Y)}
\def \jkxy {J^k(X,Y)}
\def \gkxy {G^k(X,Y)}
\def \exy {E(X,Y)}
\def \sexy{\sccenikrE(X,Y)}
\def \hn {holomorphically nondegenerate}
\def\hyp{hypersurface}
\def\prt#1{{\partial \over\partial #1}}
\def\det{{\text{\rm det}}}
\def\wob{{w\over B(z)}}
\def\co{\chi_1}
\def\po{p_0}
\def\fb {\bar f}
\def\gb {\bar g}
\def\Fb {\ov F}
\def\Gb {\ov G}
\def\Hb {\ov H}
\def\zb {\bar z}
\def\wb {\bar w}
\def \qb {\bar Q}
\def \t {\tau}
\def\z{\chi}
\def\w{\tau}
\def\Z{\zeta}
\def\phi{\varphi}
\def\eps{\varepsilon}

\def \T {\theta}
\def \Th {\Theta}
\def \L {\Lambda}
\def\b {\beta}
\def\a {\alpha}
\def\o {\omegacenik}
\def\l {\lambda}

\def \im{\text{\rm Im }}
\def \re{\text{\rm Re }}
\def \Char{\text{\rm Char }}
\def \supp{\text{\rm supp }}
\def \codim{\text{\rm codim }}
\def \Ht{\text{\rm ht }}
\def \Dt{\text{\rm dt }}
\def \hO{\widehat{\mathcal O}}
\def \cl{\text{\rm cl }}
\def \bR{\mathbb R}
\def \bS{\mathbb S}
\def \bK{\mathbb K}
\def \bD{\mathbb D}
\def \bC{\mathbb C}
\def \C{\mathbb C}
\def \N{\mathbb N}
\def \bL{\mathbb L}
\def \bZ{\mathbb Z}
\def \bN{\mathbb N}
\def \scrF{\mathcal F}
\def \scrK{\mathcal K}
\def \mc #1 {\mathcal {#1}}
\def \scrM{\mathcal M}
\def \cR{\mathcal R}
\def \scrJ{\mathcal J}
\def \scrA{\mathcal A}
\def \scrO{\mathcal O}
\def \scrV{\mathcal V}
\def \scrL{\mathcal L}
\def \scrE{\mathcal E}
\def \hol{\text{\rm hol}}
\def \aut{\text{\rm aut}}
\def \Aut{\text{\rm Aut}}
\def \J{\text{\rm Jac}}
\def\jet#1#2{J^{#1}_{#2}}
\def\gp#1{G^{#1}}
\def\gpo{\gp {2k_0}_0}
\def\emmp {\scrF(M,p;M',p')}
\def\rk{\text{\rm rk\,}}
\def\Orb{\text{\rm Orb\,}}
\def\Exp{\text{\rm Exp\,}}
\def\Span{\text{\rm span\,}}
\def\d{\partial}
\def\D{\3J}
\def\pr{{\rm pr}}

\def \CZZ {\C \dbl Z,\zeta \dbr}
\def \D{\text{\rm Der}\,}
\def \Rk{\text{\rm Rk}\,}
\def \CR{\text{\rm CR}}
\def \ima{\text{\rm im}\,}
\def \I {\mathcal I}

\def \M {\mathcal M}
\def\5#1{\mathfrak{#1}}

\newtheorem{Thm}{Theorem}[section]
\newtheorem{Cor}[Thm]{Corollary}
\newtheorem{Pro}[Thm]{Proposition}
\newtheorem{Lem}[Thm]{Lemma}
\newtheorem{Conj}[Thm]{Conjecture}

\theoremstyle{definition}\newtheorem{Def}[Thm]{Definition}

\theoremstyle{remark}
\newtheorem{Rem}[Thm]{Remark}
\newtheorem{Exa}[Thm]{Example}
\newtheorem{Exs}[Thm]{Examples}

\def\bl{\begin{Lem}}
\def\el{\end{Lem}}
\def\bp{\begin{Pro}}
\def\ep{\end{Pro}}
\def\bt{\begin{Thm}}
\def\et{\end{Thm}}
\def\bc{\begin{Cor}}
\def\ec{\end{Cor}}
\def\bd{\begin{Def}}
\def\ed{\end{Def}}
\def\br{\begin{Rem}}
\def\er{\end{Rem}}
\def\be{\begin{Exa}}
\def\ee{\end{Exa}}
\def\bpf{\begin{proof}}
\def\epf{\end{proof}}
\def\ben{\begin{enumerate}}
\def\een{\end{enumerate}}

\newcommand{\zz}{(z,\bar z)}

\maketitle

\section{Introduction}
This paper  provides an important necessary step towards solving
the local equivalence problem for hypersurfaces of finite Catlin
multitype, by giving a full classification of their polynomial
models
with nonlinear infinitesimal CR automorphisms.
Note that by
the classical Chern-Moser theory,
the only
strongly pseudoconvex  hypersurface which admits such
automorphisms  is the sphere (see
\cite{CM}, \cite{KL}).

The  Levi degenerate case has recently attracted   considerable
attention and  has led to the discovery of new types of nonlinear
symmetries (see e.g. \cite{FK}, \cite{KoLa}, \cite{KM1}).  In
contrast to the $\mathbb C^2$ case, nonlinear
infinitesimal CR automorphisms
 with coefficients of arbitrarily high degree may arise in
 $\mathbb C^n, n>2$ (\cite{ELZ2}, \cite{FM}).


Our aim in this paper  is to study systematically
nonlinear infinitesimal CR automorphisms of Levi degenerate
hypersurfaces in complex dimension three. The results provide a
description  of hypersurfaces of finite Catlin multitype  in $\Bbb
C^3$, whose models admit such automorphisms.
Since the Lie algebra of infinitesimal CR automorphisms of a
polynomial model is in one-to-one correspondence with the kernel
of the generalized Chern-Moser operator (\cite{KMZ}),
we can prove a sharp 1-jet determination result for the
biholomorphisms of any such hypersurface.

Moreover, we  identify the common source of such automorphisms.
In all cases they arise from suitable holomorphic  mappings into a
quadric in $\mathbb C^K$, $K\geq 3$,  as a ``pull-back'' of an
automorphism of the quadric.

As an application, we will determine all possible dimensions of
the Lie algebra of infinitesimal CR automorphisms for such models,
which reveals a "secondary" gap phenomenon at dimension eight.

    Our starting point are the sharp results of \cite{KMZ}
which give  an effective bound for the weighted degree of the coefficients
of an infinitesimal CR automorphism, and put restrictions on the possible form of
such vector fields.


Our first result deals with
a holomorphically nondegenerate  model hypersurface $M_H$ given by a
homogeneous polynomial $P$ of degree $d>2$ without pluriharmonic terms,
\begin{equation}\label{model1}
M_{H}:=\{\Im w = P(z,\bar z)\}, \quad (z,w)\in \C^2\times\C.
\end{equation}
Notice that the case of $d=2$ corresponds to  Levi nondegenerate
models, i.e. hyperquadrics.  We showed in  \cite{KMZ} and in
\cite{KM1} that the Lie algebra $\5g = \aut(M_H, 0)$ of all germs
of infinitesimal automorphisms of $M_H$ at $0$ admits the weighted
grading
\begin{equation}\label{sid}
\5g = \5g_{-1}\oplus \5g_{-1/d}  \oplus \5g_{0}
\oplus \5g_{\tau/d} \oplus \5g_{1-1/d} \oplus \5g_{1},
\end{equation}
 for some integer $\tau$, where $1 \le \tau \le d-2.$
The following theorem shows that if $\dim \5g_{1-1/d}\ne 0,$  then
there is a unique choice for $M_H.$

\bt
Let $M_H$ be the holomorphically nondegenerate  model hypersurface given by \eqref{model1} with $d>2$, and let
$\5g_{1-1/d}$
in
\eqref{sid} satisfy  \begin{equation}\dim \5g_{1-1/d}>0.\end{equation}
Then   $M_H$ is biholomorphically
equivalent to
\begin{equation}
\Im w =  \Re z_1 \bar z_2^l .
\label{ex1}
\end{equation}
\et

Further  we consider the more general case of
 a  holomorphically nondegenerate weighted homogeneous
model
of finite Catlin multitype.
Let
\begin{equation}\Label{model}
M_H:=\{\Im w = P_C(z,\bar z)\}, \quad (z,w)\in \C^2\times\C,
\end{equation}
where $P_C$ is a weighted homogeneous polynomial of degree one with respect to the multitype weights $\mu_1, \mu_2$
(see Section 2 for the needed definitions).

As proved in \cite{KMZ}, the   Lie algebra of infinitesimal automorphisms
$\5g=\aut(M_H,0)$ of $M_H$ admits the  weighted decomposition given
by
\begin{equation}\label{sid2}
\5g = \5g_{-1} \oplus \bigoplus_{j=1}^{2}\5g_{- \mu_j} \oplus \5g_{0}
\oplus \5g_c \oplus  \5g_{n}
 \oplus \5g_{1},
\end{equation}
where $\5g_{c}$ contains  vector fields commuting with
$W = \d_w$ and
 $\5g_{n}$ contains vector fields not commuting with $W$, whose weights in both cases
 lie
 in the interval $(0,1)$.
Note that by a result of \cite{KMZ}, vector fields in $\5g_j$ with $j<0$ are regular and vector fields
in $\5g_0$ are linear.

 \bt \label{th12}
Let $P_C(z,\bar z)$ be a weighted homogeneous polynomial of degree
1 with respect to the multitype weights,
such that the hypersurface
\begin{equation}\Label{model}
M_H:=\{\Im w = P_C(z,\bar z)\}, \quad (z,w)\in \C^2\times\C,
\end{equation} is holomorphically nondegenerate. Let
$\5g_{n}$
in \eqref{sid2} satisfy \begin{equation}\dim \5g_{n}>0.\end{equation}

Then   $M_H$ is biholomorphically
equivalent to
\begin{equation}
\Im w =  \Re z_1 \bar z_2^l
\end{equation}
or
\begin{equation}
\Im w =  \vert z_1 \vert^2 \pm \vert z_2 \vert^{2l}.
\label{ex2}
\end{equation}
\et

Note that the Levi nondegenerate case, corresponding to $l=1$, is covered by Theorem 1.2.
The following result, which deals with the
component $\5g_c$ was obtained in
\cite{KM1}.

\bd
Let $Y$ be a weighted homogeneous vector field.
A pair of finite sequences   of holomorphic  weighted homogeneous polynomials $\{U^1, \dots, U^n\}$ and  $\{V^1, \dots, V^n\}$ is called a symmetric pair of $Y-$chains if
\begin{equation}
Y(U^n)=0, \ \ Y(U^j) =c_j U^{j+1}, \ \ j=1, \dots, n-1,
\label{yuv}
\end{equation}
\begin{equation}
Y(V^n)=0, \ \ Y(V^j) =d_j V^{j+1}, \ \ j=1, \dots, n-1,
\end{equation}
where $c_j$, $d_j$ are non zero complex constants, which satisfy
\begin{equation} c_j= - \bar d_{n-j}. \label{dcj}
\end{equation}
If the two sequences are identical
we say that $\{U^1, \dots, U^n\}$ is a symmetric $Y$ - chain.
\ed

%

The following theorem shows that in general the elements of $\5g_c$  arise
from symmetric pairs of chains.

\bt \label{pivo} Let
 $M_H$  be a holomorphically nondegenerate hypersurface  given by  \eqref{model},  which admits a nontrivial
 $Y\in \5g_c$.
Then $P_C$ can be decomposed in the following way

\begin{equation}P_C=\sum_{j=1}^M T_j,\label{Barou19}
\end{equation} where each $T_j$ is given by
\begin{equation}\label{Barou20}
 T_j =\Re (\sum_{k=1}^{N_j}
  {U_j^{k}}
  {\overline {{V_j^{N_j -k +1}}}}),
  \end{equation}
where $\{{ {{U_j^{1}}, \dots, {U_j^{N_j}} }}\}$  and  $\{{ {{V_j^{1}}, \dots, {V_j^{N_j}} }}\}$ are a symmetric pair of $Y-$ chains.


Conversely, if $Y$ and $P_C$ satisfy  \eqref{yuv} --
\eqref{Barou20},
then $Y \in \5g_c$.
\et
Note that $Y$ is uniquely and explicitely determined by $P$ (see \cite{KM1}).
Hence for a given hypersurface
this result also provides a  constructive tool to determine $\5g_c$.

\bd  If $P_C$ satisfies \eqref{yuv} --
\eqref{Barou20}, the associated hypersurface  $M_H$ is called a chain hypersurface.
\ed

The description of the remaining component $\5g_1$ is a consequence of Theorem 4.7 in \cite{KMZ}
(see section 2 for the notation).
\bd We say that $P_C$ given by  \eqref{model} is balanced if it can be written as
\begin{equation}
 P_C(z, \bar z) = \sum_{|\al|_{\Lambda}= |\bar \al|_{\Lambda}
 = 1
 } A_{\al, \bar \al} z^{\al} \bar z^{\bar \al},
 \label{bala}
\end{equation}
for some nonzero pair of real numbers  $\Lambda = (\lambda_1, \lambda_2)$,
where
$$\vert \alpha\vert_{\La} := \la_1 \al_1 + \la_2 \al_2.$$
The associated hypersurface $M_H$ is called a balanced hypersurface.

\ed
Note that $P_C$ is balanced if and only if the linear vector field
$$Y=\lambda_1 z_1 \partial_{z_1} +  \lambda_2 z_2 \partial_{z_2}$$
is a complex reproducing field in the terminology of \cite{KMZ}, i.e.,
$ Y(P_C) = P_C$.

\bt The component
$\5g_1$ satisfies $\dim \5g_{1}>0 $ if and only if in suitable multitype coordinates $M_H$ is a balanced hypersurface.


\et


 As a consequence, we obtain the following result for a general hypersurface of finite Catlin multitype.
 \bt
 Let $M$ be a smooth hypersurface and $p \in M$ be a point of finite Catlin multitype with holomorphically
 nondegenerate model. If its model at $p$ is
 neither a balanced hypersurface nor a chain hypersurface,
 then its  automorphisms are determined by  the  1-jets at $p$.
 \et

Our results also confirm a conjecture about the origin of nonlinear automorphisms of Levi degenerate hypersurfaces
 formulated by the first author (see \cite{KM1}, \cite{AIM}).
 Recall that two vector fields $X_1$ and $X_2$  are  $f$-related if $f_*(X_1) = X_2.$


\bt Let $M_H$ be a holomorphically nondegenerate hypersurface
given by \eqref{model} and $Y$ be a vector field of strictly
positive weight.  Then $Y \in \aut(M_H, 0)$,  if and only if there
exists an integer $K \ge 3$ and a holomorphic mapping $f$ from a
neighbourhood of the origin in $\mathbb C^3$ into $\mathbb C^K$
which maps  $M_H$ into a Levi nondegenerate hyperquadric
$H\subseteq \mathbb C^K$ such that
 $Y$ is $f$-related with an infinitesimal CR automorphism of $H$.
\et

Let us remark that mappings of CR manifolds into hyperquadrics have been studied intensively in recent years (see e.g.
\cite{BEH}, \cite{EHZ}).
Here we ask in addition that the mapping be compatible with a symmetry of the hyperquadric.

As an application of our results, we will determine all possible
dimensions of $\aut (M_H, 0)$. Recall that the dimension  of the
symmetry group of a  Levi nondegenerate  hyperquadric in $\mathbb
C^3$ is equal to 15. In complex dimension two, the possible
dimensions of the symmetry group of a general hypersurface are
known to be $1, 2, 3, 4, 5, 8$ (\cite{KS}).

The problem of possible dimensions of the symmetry groups of a
general hypersurface in $\mathbb C^3$ was considered by Beloshapka
in \cite{Beloshapka}. He proved that for a germ of (an arbitrary)
smooth real hypersurface, which is not equivalent to a
hyperquadric, the dimension is at most 11.
The following
result shows that for polynomial models of finite Catlin multitype
the largest possible dimension is 10. Further, it demonstrates the
existence of a secondary gap phenomenon in dimension 8.

\bt
Let  $M_H$ be a holomorphically nondegenerate hypersurface given by \eqref{model}, such that $0$ is a Levi degenerate point.
Then
$$\dim \aut (M_H, 0) \leq 10.$$
If the dimension is equal to 10, then $M_H$ is equivalent to \eqref{ex1}.
If the dimension is equal to 9, then $M_H$ is equivalent to \eqref{ex2}.
Further, the dimension of $ \aut (M_H, 0)$ can attain any value from the set
$$\{ 2, 3, 4, 5, 6, 7, 9, 10\}.$$
In particular, there is no $M_H$ with $\dim  \aut (M_H, 0)   = 8$.
\et

%
%
%
%
%
%
%

The paper is organized as follows. Section 2 contains the necessary definitions used in the rest of the paper.
Section 3 deals with the $\5g_n$ component of the Lie algebra
$\5g$. Section 4 contains the proofs of the  results, up to Theorem 1.9. Section 5
deals with  possible dimensions of $\aut (M_H,0)$ and contains the
proof of Theorem 1.10.

\section{Preliminaries}
In this section we introduce notation and recall briefly some needed
definitions (for more details, see e.g. \cite{Ko1}).

Consider a  smooth hypersurface $M \subseteq \mathbb C^{3}$ and
let $p \in M $ be a  point of { finite type} $m \ge 2$  (in the
sense of Kohn and Bloom-Graham, \cite{BG}, \cite{K}).

\noindent Let $(z,w)$ be
local holomorphic coordinates centered at $p$,
where $z =(z_1, z_2)$ and  $z_j = x_j + iy_j$, $j=1,2$, and
$w=u+iv$. We assume that the hyperplane $\{ v=0 \}$ is  tangent to
$M$ at $p$, so  $M$  is described near $p$ as the graph of a uniquely
determined real valued function
\begin{equation} v = F(z_1,z_2, \bar z_1, \bar z_2,  u),
\label{vp}
\end{equation}
where $ \ dF(0) =0.$
We can assume
that (see e.g. \cite{BG})
\begin{equation}\label{fifi}
F(z_1,z_2,  \bar z_1,\bar z_2,  u)=P_m(z, \bar z) +o(u,|z|^m),
\end{equation}
where $P_m(z, \bar z)$ is a nonzero homogeneous polynomial of degree $m$ without pluriharmonic terms.

The definition of Catlin multitype involves   rational  weights associated  to the variables
$w, z_1, z_2$.
%
%
%
%
%
 The
variables $w$, $u$ and $v$ are given weight one, { reflecting} our choice
of  tangential and normal variables.
The complex tangential variables $(z_1, z_2)$  are treated
as follows.

By a  weight we understand a pair of nonnegative
 rational numbers $\La = (\la_1,
\la_2)$, where $0 \leq\la_j\leq \frac12$, and $\la_1 \ge
\la_{2}$.
Let $\La = (\la_1,
\la_2)$ be a weight, and   $\al=(\al_1, \al_2), \ $
 $\ \beta=(\beta_1,\beta_2) $   be  multiindices.
The weighted degree $\kappa$ of a monomial $$q(z, \bar z,u)=c_{\al
\beta l}z^{\al}\bar z^\beta u^{l} , \ l \in \mathbb N,$$ is then
$$ \kappa:=
l +  \sum_{i=1}^2 (\al_i + \beta_i ) \la_i.$$

%

{A polynomial $Q(z, \bar z, u)$  is weighted homogeneous
of weighted degree $\kappa$ if it is a sum of
 monomials of weighted degree $\kappa$.}

\bd

For a weight  $\La$,  the weighted length of a multiindex $\al = (\al_1, \al_2)$ is
defined by

$$\vert \alpha\vert_{\La} := \la_1 \al_1 + \la_2 \al_2.$$

Similarly, if $\al = (\al_1, \al_2)$ and  $\hat \al =
(\hat \al_1, \hat \al_2)$ are two multiindices, the weighted
length of the  pair $(\al, \hat \al)$ is
$$\vert (\alpha,\hat \al) \vert_{\La} := \la_1 (\al_1 +\hat \al_1) +
\la_2 (\al_2 + \hat \al_2).$$
\ed


\bd{ A weight $\La$ will be called distinguished for $M$ if there exist
local holomorphic coordinates $(z,w)$ in which the defining equation of $M$ takes form
\begin{equation} v = P\zz + o_{\La}(1),
\label{1}
\end{equation}
where $P\zz$ is a nonzero $\La$ - homogeneous polynomial of
weighted degree $1$ without pluriharmonic terms, and $o_{\La}(1)$
denotes a smooth function whose derivatives of weighted order less than or equal to
one vanish.}
\ed

Distinguished weights do always exist, as follows from \eqref{fifi}.
For these coordinates $(z,w),$ we have
$\Lambda=(\dfrac{1}{m}, \dfrac{1}{m}).$

In the following we shall consider the standard lexicographic order on the set
of pairs.
We recall the following definition (see \cite{C}).

\bd \label{2.6}
Let  $\Lambda_M = (\mu_1, \mu_2)$  be the infimum of all possible
distinguished weights  $\Lambda$ with respect to the lexicographic order.
The multitype of $M$ at $p$ is defined to be the pair $$(m_1,
 m_2),$$ where
$$m_j = \begin{cases}   \frac1{\mu_j} \ \  {\text{ if}} \ \  \mu_j \neq 0\\
  \infty \ \ {\text{ if}} \ \   \mu_j = 0.
\end{cases} $$
\ed
 If none of the $m_j$ is
infinity, we say that $M$ is of {finite multitype at $p$}.
Clearly, since the definition of multitype  includes all
distinguished weights, the infimum is a { biholomorphic invariant}.

{Coordinates corresponding to the multitype weight $\Lambda_M$, in
which the local description of $M$ has form (\ref{1}), with $P$
being  $\Lambda_M$-homogeneous, are called
{ multitype coordinates}.}

\bd  Let $M$ be given by (2.3).
We define a
 model hypersurface $M_H$ associated
to $M$ at $p$ by
\begin{equation} M_H = \{(z,w) \in \mathbb C^{n+1}\ | \
 v  = P_C \zz \}. \label{22}\end{equation}\ed


%
%
%


%

 Next let us  recall the following definitions.

 \bd{Let $X$ be  a holomorphic vector field
 of the form
 \begin{equation}
 X = \sum_{j=1}^2 f^j(z,w) \partial_{ z_j} + g(z,w)\partial_{w}.
 \end{equation}
 We say that  $X$  is rigid if $f^1,  f^2, g $ are all independent of the variable  $w$.}
 \ed


 We can divide homogeneous rigid  vector fields into three types, and introduce
 the following terminology.

 \bd
 Let $X \in \aut(M_H,0)$ be a rigid weighted homogeneous vector
 field. $X$ is called
 \begin{enumerate}
 \item a  shift if the weighted degree of $X$ is less than zero;
 \item a  rotation
 if the weighted degree of $X$ is equal to  zero;
 \item  a
  generalized rotation  if the weighted degree of $X$ is bigger than
  zero and less than one.
 \end{enumerate}
 \ed

 \section{Computing  $\5g_n$}

We consider a holomorphically nondegenerate model hypersurface
$M_H$  of finite Catlin multitype  in $\mathbb C^3$ , given by
\eqref{model}.
Our aim is to find all such hypersurfaces
which posses  nontrivial $\5g_n$.

For easier notation, we will write $P$ instead of $P_C$.
If $M_H$ has nontrivial $\5g_{-\mu_j}$ and $X\in \5g_{-\mu_j}$, then
by Lemma 6.1 in \cite{KMZ} there exist local holomorphic coordinates
preserving the multitype (with pluriharmonic terms allowed), such that
\begin{equation}
  X=i \d_{z_{j}}.
  \label{xi}
\end{equation}
Permuting coordinates, if necessary, we will assume that $j=1$ (hence we allow $\mu_1 < \mu_k$ for some $k$).

It follows from the assumption that we may write $P$ as

\begin{equation}
 P(z, \bar z) = \sum_{j=0}^{m} x_1^j P_j(z', \bar z'),\ \
 \label{pzz}
\end{equation}
where $P_j(z', \bar z')$ are homogeneous real valued polynomials in the variables $z' = (z_2, \dots, z_n) $  and  $ P_m \neq 0.$

\bl\Label{first-nt}
 Let $X=i \d_{z_{1}}$
be in
$ aut(M_H, 0)$ and $P$ be of the form \eqref{pzz}.
If there is a vector field $Y$ in $aut(M_H, 0)$ such that
$[Y,W]=X$, then $ m \leq 2$, i.e.  $P$ has the form
\begin{equation}\label{hopo}
 P(z, \bar z) = x_1^2P_2(z', \bar z')  + x_1 P_1(z', \bar z') +  P_0(z', \bar z')
  .
\end{equation}
\el

\bpf

  Suppose, by contradiction, that $m>2.$
We split $Y$ according to the powers of $z_1$, writing \begin{equation}\label{ho}
Y= iw\partial_{z_1}+\sum_{j=-m}^k{}Y_j,
\end{equation}
 where $Y_j$ is of the form
\begin{equation}
Y_j={\phi_1}^j(z') {z_1}^{j+1}\partial_{z_1}+ \sum_{l=2}^n {\phi_l}^j(z') {z_1}^{j}\partial_{z_l}+{\psi}^j(z') {z_1}^{j+m}\partial_{w},\end{equation}
with $ {\phi_1}^{j-1}(z')={\phi_k}^j(z')=0$ for   $j<0$
and $Y_k \neq 0$.
Let us denote
\begin{equation}
 Y' = \sum_{l=2}^n {\phi_l}^j(z') {z_1}^{j}\partial_{z_l}.
\end{equation}

Note that by weighted homogeneity of $Y$, each coefficient is homogeneous in $z'$.

We claim  that $2m-1\le m+k.$ Indeed, if not,  applying $\Re Y$ to $P-v $
the first term  of the right handside of \eqref{ho} gives
$$ - \frac{m}2 P_m^2 {x_1}^{2m-1}$$
while
all  other terms are of maximal power
$m+k$
with respect to  $z_1$, a contradiction.

We will next show that ${\psi}^k(z')=0$.
Indeed, consider the leading term with respect to the variable $z_1$ in the tangency equation
$\Re Y (P-v) = 0$.
We obtain
\begin{equation}\label{hoho2}
 \delta_{k, m-1} mx_1^{2m-1} P_{m}^2 =  mx_1^{m-1} P_m  \Re \phi_1^k z_1^{k+1}
 + 2x_1^m
Y'(P_m)
- \Im  \psi^k z_1^{k+m},
\end{equation}
where $\delta$ is the  Kronecker symbol.
We observe that, since $m>2,$  ${\Im \psi}^k(z') {z_1}^{k+m}$ cannot contain
terms in $x_1 {y_1}^{m+k-1}$ and in ${y_1}^{m+k}.$ Hence
 $\psi^k(z')=0.$

We will further consider two cases.
\newline
1. Let $Y'(P_m) \neq 0.$ We claim that  $ k = m-1$.
 Indeed, if $ k > m-1$, then we have
\begin{equation}\label{hoho33}
0 =  mP_m  \Re \phi_1^k z_1^{k+1}
 + 2x_1 Y'(P_m).
\end{equation}
Since $m>2$ and $k>1$,  the first term is harmonic in $z_1$, but the
second one is not. That gives a contradiction.
Now, using $k = m-1$, we
compare degrees in $z'$ in \eqref{hoho2}. By homogeneity it follows that $ \phi_1^k $ and  $P_m$ have the same
degree.
Dividing by $x_1^{m-1}$ and looking at coefficients of $y_1^{k+1}$ we obtain that
$\phi_1^k $ is a constant, hence $P_m$ is a constant, which gives a contradiction.
\newline
2. Let $Y'(P_m) = 0.$
We obtain
\begin{equation}
  m  x_1^{2m-1} =  m x_1^{m-1}  \Re \phi_1^{m-1} z_1^{m}
 + 2 x_1^{2m-1-k}  \Re \beta z_1^k
\frac{\partial P_{2m-1-k}}{\partial z'} -
\Im  \psi^{m-1} z_1^{2m-1}.
\end{equation}
Hence the result is a consequence of the following lemma.
\bl
There exist uniquely determined complex numbers $\al_1, \dots, \al_{m-1}$
such that
\begin{equation}
 x^{2m-1} = \sum_{j=1}^{m-1} x^{j} \Re \al_{2m-1-j} check
\end{equation}
\bpf
By comparing terms of ... we obtain the value of $\al_1$. Continuing we obtain other values of $\alpha_j$
\epf
\el

 This achieves the proof of the lemma.
\epf

\bl\Label{second-nt}
Let $X=i \d_{z_{1}}$
be in
$ aut(M_H, 0)$ and $P$ be of the form \eqref{pzz}, with $m =2$.
If there is a vector field $Y$ in $aut(M_H, 0)$ such that
$[Y,W]=X$, then $P_m$ is constant, hence $P$ has the form
\begin{equation}\label{hopo}
 P(z, \bar z) = C x_1^2  + P_1(z', \bar z') x_1   + P_0(z', \bar z')
\end{equation}
for some nonzero real constant  $C$.
\el

\bpf
Let
$k$ be given by \eqref{ho}. First assume $k=1$.
From coefficients of degree three   with respect to $z_1,$
we obtain,
\begin{equation}\label{ho2aa}
 2x_1^{3} P_{2}^2 =  2x_1 P_2  \Re \phi_1^1 z_1^2  + 2 x_1^2 \Re z_1 Y_1'(P_2)
\frac{\partial P_{2}}{\partial z'}
- \Im \psi^1 z_1^3.
\end{equation}

If $\psi^1(z') \ne 0$ in \eqref{ho2aa}, then it is a constant,
by comparing terms in $y_1^{3}.$
Hence, by homogeneity, $P_2$ is constant.


Next, assume  that  $\psi^1(z') = 0.$
Comparing degrees in $z'$, we see that $ \phi_1^1 $ and  $ P_2$ have the same
degree, or $\phi_1^1 =0. $ If $\phi_1^1 \neq 0, $then
from the coefficients of $x_1 y_1^{2}$ we obtain that
$\phi_1^1 $ is a constant, hence $P_2$ is a constant.  On the other hand,
$ \phi_1^1 =0$ implies
\begin{equation}\label{ho2}
 x_1^{3} P_{2}^2 =  x_1^2 \Re z_1 Y_1'(P_2)
\end{equation}
which is impossible, since by positivity, the left hand side contains a nonzero balanced term in $z'$, while
the right hand side has no such terms.
since  $\deg \phi_2^1 = \deg P_2 + 1$  weights??.

Now assume that $k \neq 1$.
For terms of degree $k+2$ we get
\begin{equation}\label{ho3}
0 =   2x_1 P_2  \Re \phi_1^{k} {z_1}^{k+1}  + 2 x_1^2 \Re \phi_2^k z_1^k
\frac{\partial P_{2}}{\partial z'}  -
  \Im \psi^k {z_1}^{k+2}.
\end{equation}

From the coefficients
of  $y_1^{2+k},$ we  see that $\psi^k(z')$ is a constant. If $\psi^k(z')\neq 0, $ then by homogeneity $P_2$
is a constant. Next, let $\psi^k(z')= 0.$ After dividing by $x_1$, the first term is harmonic in $z_1$, while the
second one is not, unless $k=0$. But we know from the proof of the previous lemma that $k \geq m-1$,
hence $k = 0 $ implies $P_2=0$.
This  achieves the proof of the lemma.
 \epf

The following lemma considers the case $m=1$.
Notice that for $m=1$ the definition of mutitype implies
$\mu_1= \mu_2$, dale uz nevime...

\bl\Label{third-nt}
Let $X=i \d_{z_{1}}$
be in
$ aut(M_H, 0)$ and $P$ be of the form \eqref{pzz} with $m =1$. Let $P_0$
contain no harmonic terms.
There is a vector field $Y$ in $aut(M_H, 0)$ such that
$[Y,W]=X$ if and only if
\begin{equation}
 P(z, \bar z) =  x_1 \Re \al z'^{l+1}
\end{equation}
for some  $l\in \mathbb N$, $\al \in \mathbb C$, or
\begin{equation}
 P(z, \bar z) =  x_1 \Re \al z' + \varepsilon \vert z' \vert^2,
\end{equation}
where $\epsilon \in \mathbb R$ and $\al \in \mathbb C$.
\el

\bpf
$Y$ has to be again of the form
\begin{equation}
 Y = i  w \d_{z_1} + \sum_{j=1}^2 \phi_j \d_{z_j} + \psi \d_w.
\end{equation}
From $\Re Y(P-v) = 0$, using $\Re X(P) = 0$, we obtain
\begin{equation}
 P_0 P_1 + x_1 P_1^2 = 2 x_1 \Re \phi_2 \frac{\partial P_{1}}{\partial z'} +
  \Re \phi_1
 P_1 +
2 \Re \phi_2  \frac{\partial P_{0}}{\partial z'}
-  \Im \psi.
\end{equation}
Hence for the constant and linear terms in $z_1$ we have
\begin{equation}
 P_0 P_1 = 2 \Re  \phi_2^0  \frac{\partial P_{0}}{\partial z'} +
 \Re \phi_1^{-1} P_1 -
\Im \psi^{-1}
\label{x10}
\end{equation}
and
\begin{equation}
 x_1 P_1^2 = 2 x_1 \Re  \phi_2^0  \frac{\partial P_{1}}{\partial z'} +
 \Re \phi_1^0 z_1  P_1 + 2 \Re
\phi_2^1 z_1 \frac{\partial P_{0}}{\partial z'} -
\Im \psi^0 z_1.
\label{x11}
\end{equation}
Let $k$ be as in \eqref{ho} and let
first
$k >0$.
We get

\begin{equation}
 0 = \Re \phi_1^{k} z_1^{k+1} P_1 + 2 x_1 \Re   \phi_2^{k} z_1^{k}
 \frac{\partial P_{1}}{\partial z'}
 -  \Im \psi^{k}z_1^{k+1}.
\end{equation}
From  the coefficient of $\bar z_1 z_1^{k}$, we obtain  that the middle term is zero.
It follows that
 $\Re \phi_1^{k} z_1^{k+1} P_1$ is pluriharmonic,
hence
$P_1$
is constant, which is impossible.


Now let $k=0$. We get
\begin{equation}
 x_1 P_1^2 = 2 x_1 \Re  \phi_2^0  \frac{\partial P_{1}}{\partial z'} +
 \Re \phi_1^0 z_1 P_1 - \Im  \psi^{0} z_1.
\label{x11bb}
\end{equation}
From the coefficients of $y_1$, we get
\begin{equation} - \Im \phi_1^0 P_1 - \Re  \psi^{0}= 0. \label{yy11}
\end{equation}

 This implies  that $P_1$ is harmonic, namely
$P_1 = c \Re \phi_1^0$ for some $c \in \mathbb R$.
Notice that
$ \phi_1^0 =0$ leads to contradiction. Indeed,
if $ \phi_1^0 =0$, then $\psi^0 =0$,
since $P_1$ cannot be constant.
It follows that
\begin{equation}\label{ho2}
  P_{1}^2 =  2 \Re \phi_2^0
\frac{\partial P_{1}}{\partial z'},
\end{equation}
which is impossible, since the left hand side contains a nonzero balanced term
 in $z'$,
while
the right hand side has no such terms, since $\deg \phi_2^0 = \deg P_1 + 1$. That gives the contradiction.
Next consider the equation for the coefficients of $x_1$ in \eqref{x11bb},
\begin{equation}
  P_1^2 = 2  \Re  \phi_2^0  \frac{\partial P_{1}}{\partial z'} +
 \Re \phi_1^0 P_1 - \Im  \psi^{0}.
\label{x111bbb}
\end{equation}
Substituting $P_1 = c \Re \phi_1^0$, from the mixed terms we obtain  $c=1$.
Denote  $P_1 =  \Re \phi_1^0 = \Re \al z'^l$.
Notice that
the degree of $P_0$ is $l+1$, by weighted homogeneity, since by the definition of the Catlin multitype $\mu_1 = \mu_2$.
For terms of order zero we
obtain
\begin{equation}
  P_0 P_1=   \Re \phi_1^{-1} P_1 +
2 \Re  \phi_2^0  \frac{\partial P_{0}}{\partial z'}
-
\Im \psi^{-1}.
 \label{newx11}
 \end{equation}
Using the form of $P_1$, we obtain
\begin{equation}
 P_0 \Re (\alpha z'^{l}) = \Re (\delta z'^{l+1}) \Re (\alpha z'^{l}) + \Re (\beta z'^{l+1}
  \frac{\partial P_{0}}{\partial z'}) +  \Im \gamma z'^{2l+1}\label{22}
\end{equation}
for some $\beta, \gamma, \delta \in \mathbb C$. In particular, $2 \phi^0_2 = \beta z'^{l+1}$.
We write $P_0$ as
\begin{equation}
 P_0(z', \bar z') = \sum_{j=j_0}^{l+1-j_0} A_j z'^j \bar z'^{l+1-j}
\end{equation}
with $A_{j_0} \neq 0$. Recall that $j_0 \neq 0$ by assumption.
Substituting into (\ref{22}) and comparing coefficients of $z'^{j_0} \bar z'^{2l+1-j_0}$
and $z'^{l +1 - j_0} \bar z'^{l+ j_0}$ we obtain
\begin{equation}
 \alpha = j_0 \beta, \ \ \ \ \  \al = (l+1 - j_0)  \beta.
\end{equation}
That gives
\begin{equation}
 2 j_0 = l+1
\end{equation}
and $ \beta  = \frac{2}{l+1} \al$.
It follows that $l$ is odd and $P_0 = d \ab {l+1}$ for some $ d \in \mathbb R$.
If $d \neq 0$,  we use the explicit forms of $P_0$ and $P_1$ along with
\eqref{x11bb} and \eqref{newx11}.
By \eqref{yy11} we have
$$ - \frac12 \Im (\phi_1^0)^2 = \Re \psi^0,$$
which gives
$$ \psi^0 = \frac{i}2 ( \phi^0_1)^2.$$
From \eqref{x111bbb} we obtain from the equation for the coefficient of $z'^{2l}$

$$ \frac12 \al^2 = \frac{l}{l+1} \al^2,  $$
which gives $l=1$.
The converse part of the statement is immediate to verify, using the above calculations
(see also Section 4).
That finishes the proof.


\epf

 Next we turn to the second case, $P_2 \neq 0$. Using scaling in the $z_1$ variable, we may assume
 $P_2= 1$.

\bl\Label{fourth-nt}
Let $X=i \d_{z_{1}}$
be in
$ aut(M_H, 0)$ and $P$ be of the form
\begin{equation}\label{hihi}
  P(z, \bar z) = x_1^2 + x_1 P_1(z', \bar z')  + P_0(z', \bar z').
\end{equation}
Then
there is a vector field $Y$ in $aut(M_H, 0)$ such that
$[Y,W]=X$, if and only if
 $P$ is biholomorphically equivalent, by a change of multitype coordinates, to
\begin{equation}
 P(z, \bar z) =
 x_1^2 +
 c \vert z' \vert^{2l}
\end{equation}
for some $c \in \mathbb R \setminus\{0\}$ and  $l \in \mathbb N$.
\el
\bpf
Since
$ X = i \d_{z_1},$
we have again
\begin{equation}
 Y =i  w \d_{z_1}
 + \phi_1 \d_{z_1}+ \sum_{j=1}^{n} \phi_j \d_{z_j} + \psi \d_w.
\end{equation}
Without any loss of generality, we can assume that both $P_1$ and $P_0$
contain no pluriharmonic terms.
Note that pluriharmonic terms in $P_1$ can be eliminated
by a change of variables
$z_1^* = z_1 + S(z')$, where $S$ is a holomorphic polynomial in $z'$.

In the first part of the proof, we will show that under this assumption, $P_1 = 0$.
Applying $\Re Y$ to $P-v$ gives
\begin{equation}\label{hahu}
 -(2 x_1 +  P_1 )( x_1^2 + x_1 P_1
 + P_0)
 + 2 \Re \phi_1 \frac{\partial P}{\partial z_1} + 2 \Re
\sum_{j=2}^n \phi_j \frac{\partial P}{\partial z_j} -  \Im \psi = 0.
\end{equation}
Let $k$ be as in \eqref{ho}. Assume first that $k>1$.
If $\phi_1^{k}\neq 0$ we get
$$ x_1 \Re\phi_1^{k}z_1^{k+1} - \frac12 \Im \psi^k  z_1^{k+2} = 0,$$
 which gives a contradiction, since the second term is harmonic in $z_1$, while the first one is not.
Hence $\phi_1^{k}= 0$ and  we have, for terms of degree $k+1$ in $z_1$,

\begin{equation}
  2\delta_{k,2} x_1^3  = 2x_1 \Re \phi_1^{k-1} z_1^{k} + 2 x_1 \Re  \sum_{j=1}^n  \phi_j^{k} z_1^{k}
 \frac{\partial P_{1}}{\partial z_j}
 -  \Im \psi^{k-1}z_1^{k+1}.
\end{equation}
Looking at the coefficients of $y_1^{k+1}$ we see that $\psi^{k-1}$ is a real or imaginary
constant,
depending on the parity of $k$.
If
$\psi^{k-1}  \neq 0$,
then by homogeneity $P_1$ has degree at most one, hence it must be zero.
If $\psi^{k-1}=0$, after dividing by $x_1$, the right hand side is pluriharmonic in $z_1$. It follows that
$\delta_{k,2}= 0, $ hence $k > 2$.
But then $\sum_{j=1}^n  \phi_j^{k} z_1^{k}
 \frac{\partial P_{1}}{\partial z_j}$ cannot contain any mixed terms,
hence it has to vanish.
So
$$ Y'_k(p_1) = 0.$$
From the next equation we obtain
\begin{equation}
  2\delta_{k-1,2} x_1^3  = 2x_1 \Re \phi_1^{k-2} z_1^{k-1} + 2  \Re  \sum_{j=1}^n  \phi_j^{k} z_1^{k}
 \frac{\partial P_{0}}{\partial z_j} +  2  x_1\Re  \sum_{j=1}^n  \phi_j^{k-1} z_1^{k-1}
 \frac{\partial P_{1}}{\partial z_j}
 -  \Im \psi^{k-1}z_1^{k+1}.
\end{equation}
Looking at coefficients of $y_1$, we obtain that
$ 2  \Re  \sum_{j=1}^n  \phi_j^{k} z_1^{k}
 \frac{\partial P_{0}}{\partial z_j}$ is pluriharmonic, hence it has to vanish.
 That contracicts holomorphic nondegeneracy of $M_H$.
It follows that
as claimed.
\\
Now let us assume that $k=1$. For the third order terms in $z_1$ we obtain
\begin{equation}
2 x_1^3  = 2x_1 \Re \phi_1^{1} z_1^{2}
 - \Im \psi^{1}z_1^{3}.
\end{equation}
This determines $\phi_1^{1}$
and $\psi^{1}$, which are thus constant.
Looking at terms of second order in $z_1$ we obtain from \eqref{hahu}
and \eqref{hihi}
\begin{equation}
 - 3x_1^2 P_1
 + 2x_1  \Re \phi_1^0 z_1
 + 2x_1 \Re \sum_{j=2}^n \phi_j^1 z_1 \frac{\partial P_{1}}{\partial z_j} + 2\Re \phi^1_1 z_1^2 P_1 -
 \Im \psi^{0}z_1^{2} =0.
\end{equation}
Looking at coefficients of $y_1^2$, we obtain that $P_1$ is pluriharmonic, hence $P_1 =0$.

We will further assume that $P_1=0$.
We obtain  for terms linear in $z_1$,



\begin{equation}
- 2 x_1 P_0  + 2 \Re  z_1  \sum_{j=2}^n
\phi_2^1 \frac{\partial P_0}{\partial z_j} + 2x_1 \Re \phi_1^{-1} -  \Im \psi^{-1} z_1= 0
\end{equation}
which gives equations for coefficients of $x_1$ and $y_1$.  Namely
\begin{equation}
 2 P_0  + 2 \Re  \sum_{j=2}^n
\phi_j^1 \frac{\partial P_0}{\partial z_j} + 2  \Re \phi_1^{-1} -  \Im \psi^{-1}= 0
\end{equation}
and
\begin{equation}
   - 2\Im  \sum_{j=2}^n
\phi_j^1 \frac{\partial P_0}{\partial z_j} -  \Re \psi^{-1}= 0.
\end{equation}
Using $\phi_j^1 = \alpha z'$ for some $\alpha \in \mathbb C$ and the fact that $P_0$ contains no harmonic terms,  it follows that

\begin{equation}
   \Im  \sum_{j=2}^n
\phi_j^1 \frac{\partial P_0}{\partial z'} = 0,
\end{equation}
hence
\begin{equation}
 P_0 =  \sum_{j=2}^n
- \phi_j^1 \frac{\partial P_0}{\partial z'}.
\end{equation}
It follows that $P_0$ has a complex reproducing field, hence $P_0$
is a balanced polynomial.
as claimed.
That finishes the proof.
\epf

\section{Proofs of Theorems 1.1 - 1.9.
}
In this section we complete the proofs of the results stated in the introduction, up to Theorem 1.9.
Theorem 1.1  is an immediate consequence of Theorem 1.2.

{\it Proof of Theorem 1.2.}
We apply Lemma \ref{first-nt} - \ref{fourth-nt}.
Note that by  Lemma \ref{first-nt} and Lemma \ref{second-nt}, we obtain either $\mu_1 = \frac12$, or $\mu_1 = \mu_2$.
Hence we can assume without
any loss of generality that $i \d_{z_1} \in \aut(M_H, 0)$.
Using suitable scaling,
rotation and adding harmonic terms leads to the form given in the statement of the theorem. \qed
 \\[2mm]

Theorem 1.7 follows immediately from Theorem 4.7 in \cite{KMZ}.
Combining Theorem 1.2, 1.4 and 1.7 with the results of \cite{KMZ} leads to Theorem 1.8.

{\it Proof of Theorem 1.9.}
If $\5g_1 \neq 0$, then by Theorem 4.7 in \cite{KMZ} we have


\begin{equation}
 P(z, \bar z) = \sum_{|\al|_{\Lambda}= |\bar \al|_{\Lambda} = \frac12} A_{\al, \bar \al} z^{\al} \bar z^{\bar \al}
\end{equation}
for some pair $\Lambda = (\la_1, \la_2)$ (not necessarily equal to the multitype weight).
Let $K$ be the number of nonzero terms in the sum. We order the multiindices and write
$P$ as
\begin{equation}
 P(z, \bar z) = \sum_{j=1}^{K}  A_j
 z^{\al_j} \bar z^{\al_{K+j}}
\end{equation}
Consider the hyperquadric in  $\mathbb C^{2K+1}$ defined by
\begin{equation}
 \Im \eta = \sum_{j=1}^{K}
  A_j \zeta_{j}
  \overline{\zeta_{{K+j}}},
\end{equation}
and the mappping  $f : \mathbb C^3 \to \mathbb C^{2K+1}$
given by $\eta = w$ and
 $\zeta_{j} = z^{\al_j}$ for $ j = 1, \dots, 2K$.

It is immediate to verify that the vector field in $aut(M_H, 0)$
$$Y =\left( \lambda_1 z_1  \d_{z_1}+ \lambda_2 z_2  \d_{z_2} \right) w+
\dfrac{1}{2}w^2 \d_{w},\
$$
is $f$-related to the infinitesimal automorphism of the above hyperquadric
 given by
$$Z  =\frac12 \eta \sum_{j=1}^{2K} \zeta_j  \d_{\zeta_j}
+ \dfrac{1}{2}\eta^2 \d_{ \eta}.\ $$

Next, if $\5g_n\neq 0$ we consider the two cases from Theorem 1.2. We have $K=3$ in both cases. In the first case,
we define $f$ by $\eta = w, \ \zeta_1 = z_1, \zeta_2 = z_2^l$. We verify that the vector field

$$Y_1=aw  \d_{z_1}-i{\bar a} z_1{z_2}^l  \d_{z_1} -i
\bar a\dfrac{1}{l} {z_2}^{l+1}  \d_{z_2} + 2i \bar a {z_2}^l
w   \d_{w}, \ a\in \mathbb C.$$
in $aut(M_H,0)$  is $f$-related to the infinitesimal automorphism of the hyperquadric,
 \begin{equation}
 \Im \eta = \Re
  \zeta_{1}
  \overline{\zeta_{{2}}},
\end{equation}
 given by
  $$Z_1 =a\eta \d_{\zeta_1}-i{\bar a} \zeta_1{\zeta_2} \d_{\zeta_1} -i
\bar a  {\zeta_2^2}\d_{ \zeta_2} + 2i \bar a {\zeta_2}
 \eta  \d_{ \eta}, \ a\in \mathbb C.$$

The second case is completely analogous.
The case of $\5g_c$ follows from  Theorem 1.2  in \cite{KM1}.
This finishes the proof. \qed
\\[1mm]

\section{Dimension of $\aut(M_H, 0)$.}

In this section we will again assume that $M_H$ is a
holomorphically nondegenerate model given by \eqref{model}.   We
will first prove two auxiliary lemmata. Let us denote $\5g_t =
\5g_{- \mu_1} \oplus \5g_{- \mu_2}$,  the part of $\aut(M_H,0)$
containing complex tangential shifts.

\bl Let there exist
two regular vector fields in
$\5g_t$
 whose values at $0$ are linearly independent over $\mathbb R$,
but dependent over $\mathbb C$. Then $M_H$ is biholomorphic to
\begin{equation}
  \Im w  = C x_1 ^2 + x_1 \Re \al z_2^l + Q(z_2, \bar z_2)
\end{equation}
for some $C \in \mathbb R$, $\al \in \mathbb C$ and homogeneous
polynomial $Q$ without harmonic terms.
\el
\bpf
Let $Z_1$, $Z_2$ be such vector fields. Without loss of generality, we may assume that
$Z_1 = i \d_{z_1}$ and
$$Z_2 = \d_{z_1} + \psi(z_1, z_2) \d_w.$$
Note that this form is attained using transformations preserving
multitype coordinates (with pluriharmonic terms allowed). The commutator of $Z_1$, $Z_2$ either
vanishes, or lies in $\5g_{-1}$, i.e. it  is a real multiple of
$\d_w$. This leads to $ \psi(z_1, z_2) = C z_1 + \alpha z_2^l \ $
for some $\al \in \mathbb C$, $C \in \mathbb R$  and $l \in
\mathbb N$. From $\Re Z_1 (P) = 0$ it follows that
\begin{equation}
 P(z_1, z_2, \bar z_1, \bar z_2) = \sum_{j=0}^m x_1^j P_j(z_2, \bar z_2).
\end{equation}
From $\Re Z_2 (P) = 0$ we obtain that  $\Re  \dfrac{\partial
P}{\partial z_1}$ is pluriharmonic. It follows that
\begin{equation}
  P(z_1, z_2, \bar z_1, \bar z_2)  = C x_1 ^2 + x_1 \Re \al z_2^l  + P_0(z_2, \bar
  z_2),
\end{equation}
which finishes the proof.
\epf

Recall that by the results of \cite{KMZ} and \cite{KM1}, $\dim
\5g_{c} \leq 1$ and $\dim \5g_{1} \leq 1$. The following lemma
considers the case when both components are nontrivial.
\bl Assume that $\dim \5g_{c}= \dim \5g_{1} = 1$. Then  $\dim
\5g_0 = 3$.
\el \bpf Let $ Z \in \5g_0$ be a rotation and $Y \in \5g_{1}$ be
nonzero vector field.  By Theorem 4.7 in \cite{KMZ}, $Y$ has the form
 \begin{equation}\Label{ex5a'}
Y =  \sum_{j=1}^2 \phi_j(z) w \partial_{z_j}
 + \frac12 w^2  \partial_{w}
\end{equation}
where the $\phi_{j}$ have the complex reproducing field property
\begin{equation}\Label{ex8a'}
2\sum_{j=1}^2 \phi_j(z) P_{z_j} = P(z, \bar z).
 \end{equation}
It is immediate to verify that in Jordan normal form the linear vector field
$\sum_{j=1}^2 \phi_j(z) \d_{z_j}$
is diagonal with real coefficients.
We can thus consider multitype coordinates in which $Y$ is a real multiple of
$$ \left( \lambda_1 z_1  \d_{z_1}+ \lambda_2 z_2  \d_{z_2} \right) w+
\dfrac{1}{2}w^2 \d_{w}.\
$$
We claim that in such coordinates, $Z$ is also
diagonal. Indeed, let $X \in \5g_{c}$ 
be a
nonzero vector field. The commutator of $X$ and $Y$ is of weight bigger than
one, so $[X,Y] = 0$, by the results of \cite{KMZ}. It follows that
the pair $(\la_1, \la_2)$
is linearly independent
with the multitype weights $(\mu_1, \mu_2)$. If $\mu_1 \neq
\mu_2$, any rotation is diagonal. So we may assume $\mu_1 =
\mu_2$, which implies $\la_1 \neq \la_2$.
 The commutator of $Z$ with $Y$ has to be a real multiple of $Y$.
 Computing the commutator, it follows immediately
that $Z$ is diagonal. Next, the rotations with real coefficients
have dimension  one, the coefficients being given by $\l_1 -
\mu_1$ and $\la_2 -\mu_2$. Let us write $P$ in the form
$$ P(z_1, z_2, \bar z_1, \bar z_2) = \sum_{\vert \al, \hat \al \vert_{\Lambda_M}=1} A_{\al. \hat \al} z^{\al}\bar z^{\hat
\al}.$$
 If there exist in
addition two linearly independent imaginary rotations, we have
$\al_1 = \hat \al_1$, $\al_2 = \hat \al_2$ whenever $A_{\al. \hat
\al} \neq 0$
From the real rotation and weighted homogeneity we obtain
a unique solution for $\al_1, \al_2, \hat \al_1, \hat \al_2$. That
contradicts holomorphic nondegeneracy of $M_H$.
On the other hand, there is an imaginary rotation with
coefficients $i \la_1, i \la_2$. Hence there are two linearly
independent rotations, and therefore $\dim \5g_0 = 3. $
\epf

\bl There exist multitype coordinates in which every rotation is
linear. Moreover, $$\dim \5g_0 \leq 5.$$ \el

\bpf The first part of the statement is proved in Proposition 3.9
in \cite{KMZ}. In the normalization of Proposition 2.6 in
\cite{KMZ},  it is immediate to check that  the nonzero vector
fields given by
\begin{equation}
 Z = (a z_1 + \beta z_2) \d_{z_1} +  c  z_2 \d_{z_2},
\end{equation}
where $a, c \in \mathbb R, \beta \in \mathbb C$ do not belong to
$\5g_0$.  Hence the space of rotations is at most four
dimensional, and  $\dim \5g_0 \leq 5$.

 \epf

\noindent {\it Proof of Theorem 1.10}. \ First we consider the
case when $\dim \5g_n  > 0.$ By Theorem \ref{th12} $M_H$ is given
by \eqref{ex1} or \eqref{ex2}. The vector fields in $\aut(M_H, 0)$
of \eqref{ex1} are described explicitly in \cite{FM}, which shows
that the dimension is equal to 10. A completely analogous
computation gives all vector fields in $\aut(M_H, 0)$ for
\eqref{ex2}, and shows that in this case, $\dim \aut(M_H, 0) = 9.$
Indeed, we have $\dim \5g_{-\frac12} = 2$ and $\dim \5g_{\frac12}
= 2$. Further, since $\mu_1 \neq \mu_2$, all rotations have to be
diagonal, which gives immediately $\dim \5g_{0} = 3$. By the
results of \cite{KM1}, we have $\dim \5g_{c} = 0$. Since $\dim
\5g_{-1} = \dim \5g_{1}= 1$, we obtain $\dim \aut(M_H, 0) = 9$.

Next we consider the case when $\dim \5g_n  = 0$ and show that
in this case $$\dim \aut(M_H, 0) \leq 7.$$
\newline
We will consider all possible dimensions of $\5g_t$, namely
$0,1,2,3$. Note that  if $\dim \5g_t= 4$, by Lemma 5.1 it follows
that $\mu_1 = \mu_2 =\frac12$, hence $M_H$ is a Levi nondegenerate
hyperquadric.

Let us first assume that $\dim \5g_t = 0$. Then by Theorem 1.2,
Lemma 5.2 and 5.3 we obtain immediately $\dim \aut(M_H, 0) \leq
7$.

Next assume that $\dim \5g_t = 1$. If $\dim \5g_{c}+ \dim \5g_{1}
= 2$,  then by Lemma 5.2 $ \dim \5g_{0} \leq 3$, and since $\dim
\5g_n = 0$, this  leads to $\dim \aut(M_H, 0) \leq 7$. Now let
$\dim \5g_1 + \dim \5g_c \leq 1$. Without any loss of generality,
we can assume $Y = \d_{z_1} \in \5g_t.$
 By the partial normalization (Proposition 2.6 in \cite{KMZ}),
the five dimensional linear subspace of vector fields of the form
\begin{equation}
 Z = (\alpha z_1 + \beta z_2) \d_{z_1} +  d  z_2 \d_{z_2},
\end{equation}
where $\al, \beta \in \mathbb C, d \in \mathbb R$ has trivial
intersection with $\5g_0$. This follows from Lemma 5.3 and the
fact that if $Z \in \aut(M_H, 0)$, the commutator of $Z$ and $Y$ has to be a real multiple
of $Y$.  It implies $\dim \5g_0 \leq 4$.
Since $\dim \5g_n = 0$, we obtain
 $\dim \aut(M_H, 0) \leq 7$.

Next, let $\dim \5g_t = 2$ and
the values of two generators of  $\5g_t$ at $0$
are complex independent. We claim that $\dim \5g_0 \leq 2$.
Indeed, let
$Z_1 \in \5g_{- \mu_1}$, $Z_2 \in \5g_{- \mu_2}$ be two vector
fields, whose values at $0$ are complex independent. Since their
commutator is of weight $-\mu_1 - \mu_2 > - 1$, it must vanish. We
may therefore assume $Z_1 =i  \d_{z_1}$, $Z_2 = i \d_{z_2}$. Hence $P$
is a function of $x_1, x_2$.  Again, this form is attained by
a transformation preserving multitype coordinates (with
pluriharmonic terms allowed). Computing commutators of $Z_1$ and
$Z_2$ with a rigid element of $\5g_0$ of the form
\begin{equation} Y = (az_1 + b z_2) \d_{z_1} + (cz_1 + d z_2)
\d_{z_2} + \psi(z_1, z_2) \d_w,
\end{equation}
we see that the coefficients $a, b, c, d$ must be real.
It follows that
$$\Re Y(P-v) = (ax_1 + b x_2)  \dfrac{\partial P}{\partial x_1} + (cx_1 + d x_2) \dfrac{\partial P}{\partial x_2}
- \frac12 \Im \psi(z_1, z_2)= 0.
$$
The last term is pluriharmonic and the remaining part is a function of $x_1, x_2$. It follows that the pluriharmonic term vanishes,
and we obtain
$$ \frac{ \dfrac{\partial P}{\partial x_2}}{\dfrac{\partial P}{\partial x_1}} =
- \frac{ax_1 + b x_2}{cx_1 + d x_2}.$$ It follows that $a,b,c,d$
are determined uniquely by $P$, up to a real mutiple. Note that
the left hand side is nonconstant,  otherwise $M_H$ is
holomorphically degenerate. Hence $\dim \5g_0 \leq 2$. It follows
that $\dim \5g \leq 7$.

If $\5g_t$ contains    two complex dependent at $0$ vector fields, we use Lemma 5.1
and consider a rotation
\begin{equation}
 Y = (az_1 + b z_2) \d_{z_1} + (cz_1 + d z_2) \d_{z_2} + \psi(z_1, z_2)
 \d_w.
\end{equation}

 Let first $C\neq 0$. After scaling in $z_1$ and absorbing the
mixed term into $x_1^2$, we may assume
$$P(z_1, z_2, \bar z_1,
\bar z_2) = x_1^2 + Q(z_2, \bar z_2), $$
 where $Q$ is different from $\vert z_2 \vert^{l} $, since $\dim \5g_n = 0$.
 We obtain
\begin{equation}
 2x_1 \Re (az_1 + b z_2) + \Re (cz_1 + d z_2) \dfrac{\partial Q}{\partial z_2} - \frac12 \Im \psi(z_1, z_2)  = 0.
\end{equation}
From the coefficient of $y_1$ we obtain $\Im (c \dfrac{\partial
Q}{\partial z_2})= 0$, which implies $i c \d_{z_2}
\in \5g_t$. It follows that
$c=0$.
Hence
\begin{equation}
 2x_1 \Re (az_1) + 2x_1 \Re (b z_2) + \Re  d z_2 \dfrac{\partial Q}{\partial z_2}  - \frac12 \Im \psi(z_1, z_2)  = 0.
\end{equation}
The first three terms have different powers of $z_1$, so they
must all be pluriharmonic. It follows that $a$ is purely imaginary
and $b = d =0$. Hence $\dim \5g_0 \leq 2$ which implies
$\dim \5g \leq 6$.

Let now $C = 0$, i.e. $$P(z_1, z_2, \bar z_1, \bar z_2) =  x_1 \Re
\al z_2^l + Q(z_2, \bar z_2).$$ Without any loss of generality, we
will assume that $Q$ contains no harmonic terms, and no $(1,l)$
terms, which can be absorbed into  the first term by a change of
variables $z^*_1 = z_1 + \gamma z_2.$  We have
\begin{equation*}
 2\Re (az_1 + b z_2) \Re\al  z_2^{l} + 2l x_1 \Re(\al  z_2^{l-1} (cz_1 + d z_2)) +
 2\Re \left[ (cz_1+ d z_2) \dfrac{\partial Q}{\partial z_2} \right] - \Im \psi(z_1, z_2) = 0.
\end{equation*}
From the equation for the coefficient of $z_1 \bar z_1 z_2^{l-1}$
we obtain that $ c  = 0 $. From terms of order zero in $z_1$,
looking at the coefficient of $\bar z_2 z_2^l$ we obtain $b=0$,
since $Q$ contains no $(1,l)$ terms.
 From the coefficient of $ z_1 \bar z_2^l$
we obtain $ a  + l \bar  d = 0$. This gives $\dim \5g_0 \leq 3$
and  $\dim \5g \leq 7$.
\newline
If $\dim \5g_t= 3$, we use Lemma 5.1 to conclude that $P$ must be
biholomorphic to
\begin{equation}
  \Im w = \vert z_1 \vert^2 + (\Re z_2)^{l},
\end{equation}
where $l >2$. Since $\mu_1 \neq \mu_2$, the rotations have to be
diagonal. It is immediate to verify that the only rotations are
real multiples of $i z_1 \d_{z_1}$. Hence $\dim \5g_0 =2 $ and
 $\dim \5g \leq 7$.

To finish the proof, we give examples of models with $\aut(M_H, 0)
= 7,6,5,4,3,2.$
The hypersurface
given by
\begin{equation}
 \Im w = (\sum_{j=1}^2 \vert z_j \vert^2)^2
\end{equation}
admits the same rotations as the sphere, hence $\dim \5g_0 =5$ and
 $\dim \aut(M_H, 0)= 7$.
\newline
 Six dimensional
$\aut(M_H, 0)$ occurs  for
\begin{equation}
  \Im w = \vert z_1 \vert^2 + (\Re z_2)^{3},
\end{equation}
which has $\dim \5g_t=3 $ and  $\dim \5g_0 = 2$. \newline
Dimension five occurs for
  $$\Im w = \sum_{j=1}^2 \vert z_j \vert^4,$$
where $\dim \5g_0 = 3$.
\newline
Dimension  four occurs for
\begin{equation}
\Im w =
(\Re z_1)^3 + (\Re z_2)^4,
\end{equation}
where $\dim \5g_t=2.$
 \newline Dimension
three occurs for
\begin{equation}
  \Im w = \Re  z_1 \bar z_1^3 + (\Re z_2)^{3},
\end{equation}
with $\dim \5g_t=1.$
\newline Dimension two occurs for a generic model, e.g.,
\begin{equation}
  \Im w = \Re (\sum_{j=1}^2  z_j \bar z_j^3)
\end{equation}
This finishes the proof of the theorem. \qed

\end{document}